\documentclass[12pt]{iopart}

\usepackage{lineno,hyperref}
\usepackage{bm}
\usepackage{amssymb}
\usepackage{amsthm}
\usepackage{iopams}
\usepackage{graphicx}
\usepackage[dvipsnames]{xcolor}
\usepackage{siunitx}
\usepackage{algorithm}
\usepackage{algpseudocode}
\usepackage{hyperref}
\usepackage{calrsfs}
\DeclareMathAlphabet{\pazocal}{OMS}{zplm}{m}{n}

\newtheorem{theorem}{Theorem}[section]
\newtheorem{lemma}[theorem]{Lemma}

\modulolinenumbers[5]

\definecolor{editColor}{RGB}{0, 102, 204}

\begin{document}


\title[On Bias and Its Reduction via Standardization]{On Bias and Its Reduction via Standardization in Discretized  Electromagnetic Source Localization Problems}

\author[inst]{Joonas Lahtinen}
\address{Mathematics \& Computing Sciences, Tampere University, Tampere 33720, Finland}
\ead{joonas.j.lahtinen@tuni.fi}

\begin{abstract}
In electromagnetic source localization problems stemming from linearized Poisson-type equation, the aim is to locate the sources within a domain that produce given measurements on the boundary. In this type of problem, biasing of the solution is one of the main causes of mislocalization. A technique called standardization was developed to reduce biasing. However, the lack of a mathematical foundation for this method can cause difficulties in its application and confusion regarding the reliability of solutions. Here, we give a rigorous and generalized treatment for the technique using the Bayesian framework to shed light on the technique's abilities and limitations. In addition, we take a look at the noise robustness of the method, which is widely reported in numerical studies. The paper starts by giving a gentle introduction to the problem and its bias and works its way toward standardization.
\end{abstract}

{\bf An article corresponding to this preprint has been published in Inverse Problems 40 095002 on 22 July 2024. DOI: \href{https://doi.org/10.1088/1361-6420/ad5f53}{10.1088/1361-6420/ad5f53}}

\submitto{\IP}
\maketitle



\section{Introduction}
The observation model of the general electromagnetic source localization problem is based on the following form of the Poisson equation with Neumann boundary conditions:
\begin{eqnarray}\label{eq:Possion1}
        & \nabla\cdot (M\nabla u)=f\quad &\textnormal{in } \Omega,\:\:\:\,\\
        & {\bf n}\cdot (M\nabla u)=g\quad &\textnormal{on } \partial\Omega, \label{eq:Possion2}
\end{eqnarray}
for any $f\in L^2(\Omega)$ and $g\in L^2(\partial\Omega)$, where $M$ is the "conductivity" tensor of the domain and ${\bf n}$ denotes an outward unit normal vector of $\Omega$.
The goal is to localize the few small compact supports of the source term $f$ when $u$ is measured at a number of boundary points. Here we denote a vector of $m$ boundary measurements as ${\bf y}\in\mathbb{R}^m$. Such a problem can arise, e.g., in the neuroelectromagnetic inverse source problem \cite{Grech2008,Calvetti2009}, electrocardiography inverse problem \cite{GuofaShou2007StEI}, and inverse problems of geophysics \cite{BackusGeophysics1967} such as electrical resistance tomography \cite{JohnsonGeoPhys2015,ShiWenyangGeoPhys2023}. 
In such problems, the geometry of the domain is typically complex and inhomogeneous by its parameters, which is why the problem is linearized to obtain a numerical solution. For example, in the quasi-static neuroelectromagnetic inverse source localization problem \cite{Knosche2022}, using a realistic model of the head and folded brain structure helps to locate the brain activity \cite{VanrumsteBart2002realvsspherical,VattaFederica2010realvsSpherical}. The geometrically detailed model also helps to estimate cardiac electrical activity from voltage distributions measured on the body surface in electrocardiography inverse problem \cite{WangDafang2011}. Moreover,  in this ill-conditioned problem, the measurements obtained from the domain's boundary are typically noisy. Thus, one need to solve such a ${\bf x}\in\mathbb{R}^n$ from the linear observation model
\begin{eqnarray}
    {\bf y}=L{\bf x}+\underline{{\bf q}},
\end{eqnarray}
so that it has a set of non-zero values within the support of $f$. Here $L\in\mathbb{R}^{m\times n}$ is the system matrix, and the noise vector $\underline{{\bf q}}\in\mathbb{R}^m$ is assumed to follow zero-mean Gaussian distribution $\pazocal{N}({\bf 0},C)$ with measurement noise covariance matrix $C\in\mathbb{R}^{m\times m}$. In this paper, random variables are denoted by underlined symbols, and when we treat them as variables for optimization or use their realizations, we write them without underlining.

The support of $f$ is often assumed to consist of a few distinct points called {\em true sources}. One possible approach is to use a sparse solver to reconstruct a nearly point-like source. On the other hand, more complex algorithms pose limitations to the forward model's resolution. Thus, there is a trade-off between spatial accuracy and the sparsity of the solution, e.g, in the case of the well-known LASSO when it is implemented using Least Angle Regression algorithm \cite{Efron2004LASSOcomplexityN3rd}. Consequently, the practice has been to use computationally less intensive inverse algorithms and choose the component of the reconstructed parameter vector with the largest magnitude as the indicator of source location \cite{Manfred1999,SekiharaNagarajan2008,Grech2008,PascualMarqui2002}. This localization practice has been demonstrated to be robust and reliable with a range of inversion methods \cite{Neugebauer2022,Lahtinen2023}.

As the forward modeling incorporates a probabilistic interpretation of the noise, Bayesian modeling can be considered appropriate for regularization \cite{SomersaloKaipio2005_Bayes}. Consider a common Gaussian process regression model, also called Bayesian minimum norm estimate (BMNE) \cite{HamalainenMNE}, with a posterior composed of a Gaussian likelihood $\underline{\bf y}\mid{\bf x}\sim \pazocal{N}(L{\bf x},C)$ and a Gaussian prior $\underline{\bf x}\sim\pazocal{N}({\bf 0},\Gamma)$, where and $\Gamma\in\mathbb{R}^{n\times n}$ is a diagonal prior covariance. An estimate for the realization of ${\bf x}$ can be obtained as the {\em maximum a posteriori} (MAP)
\begin{eqnarray}\label{eq:MNE_MAP_estimate}
    \hat{\bf x}=\mathrm{arg\: max}_{\bf x\in\mathbb{R}^n} \left\lbrace \pazocal{N}(L{\bf x},C)\times \pazocal{N}({\bf 0},\Gamma)\right\rbrace,
\end{eqnarray}
which constitutes the following $L2$--$L2$ minimization problem, whose objective function is composed of a $L2$ norm-based data fit and a regularization term that can be associated with Tikhonov regularization:
\begin{eqnarray}\label{eq:MNE_minimization_problem}
    \min_{\bf x\in\mathbb{R}^n}\left\lbrace ({\bf y}-L{\bf x})^\mathrm{T}C^{-1}({\bf y}-L{\bf x})+{\bf x}^\mathrm{T}\Gamma^{-1}{\bf x}\right\rbrace.
\end{eqnarray}
This model, however, poses some challenges: Due to the maximum principle of elliptic partial differential equations \cite{ProtterWeinberger1967,JostJurgenPDE2007}, the maximum of the forward problem lies on the boundary. Consequently, it biases the estimated source locations to the part of the boundary where the measurement sensors are located. The most intuitive approach to tackle bias is to use solution weighting that utilizes some norm of the system matrix columns \cite{MNE_cardio_biased,Lin2006a,Attal2013_wMNEsLORdSPM}. Similar weighting has also been developed for a Bayesian method. In this case, the sensitivity weighting obtained through hyperparameters has been obtained by equating the proportion of the pure signal from the signal-to-noise ratio to the expected value of the power of the observed signal given by the Bayesian model when the measurements are assumed noiseless \cite{Calvetti2019AutommaticDepthWeighting}. The weighting of this type, however, could reduce the bias only depth-wise away from the measurement sensors. More comprehensive unbiasing can be achieved by weighting the regularization so that the point spread of each covered parameter is taken into account \cite{Pascual-Marqui2011_eLORETA,ElvetunNielsen2023}. In this case, complete unbiasing can be achieved with noiseless data and as the regularization parameter tends to zero. As the third known way to counter the bias, Pascual-Marqui has introduced the technique called {\em standardization} \cite{PascualMarqui2002}, where the components of an inverse solution are post-hoc weighted by the square roots of the diagonal elements of the resolution matrix, that is
\begin{eqnarray}
    R=\Gamma L^\mathrm{T}\left(L\Gamma L^\mathrm{T}+C\right)^{-1}L
\end{eqnarray}
as derived in \cite{AnMeijian2012ResMat}, to form unbiased estimations for covering the location of a single source. The resolution matrix is also known as the Backus-Gilbert resolution kernel \cite{BackusGeorge1968ResMat,GravedePeralta1997ResMat}. It defines a relation between the inverse solution and forward model parameters under linear inversion, essentially in $L2$-$L2$ minimization. It is mainly applied to neuroelectromagnetic source imaging \cite{WanXiaohong2008EleSourImgResMat} and inversion problems of geophysics \cite{AnMeijian2012ResMat}. The matrix can be used to estimate the resolution length of the estimation when the forward and inverse models are known, and its diagonal elements can be used to evaluate the resolvability of the target source \cite{AnMeijian2023ResMat}. Here we define localization bias as in \cite{SekiharaNagarajan2008}, by saying that the solution is unbiased when the point-spread function has its maximum at the correct source location. This means that the $k$th resolution matrix column has its maximum component at $k$th entry, as we will see later. This definition assumes no uncertainty related to the forward model $L$ as the approximation of the Poisson equation (\ref{eq:Possion1}) with the boundary conditions (\ref{eq:Possion2}). All the uncertainty comes from the random variables $\underline{\bf x}$ and $\underline{\bf q}$.

Due to the formulation of the weights, the standardization weighting procedure is closely related to the standardization of Gaussian random variables as the marginal distribution of MAP for the $k$th entry is $\pazocal{N}({\bf 0},R_{kk})$. The localization of the source via standardization of Gaussian process regression has been previously shown to be simultaneously unbiased and errorless under specific conditions \cite{Pascual2007sqrtm,SekiharaKensuke2005Lbas}, e.g., using the boundary of a volumetric structure as the domain for calculations. By errorless, we mean that the entry with the maximum magnitude of MAP estimate corresponds to the location of the true source. Moreover, the robustness of the method to the measurement noise was described in \cite{DümpelmannMatthias2012sard,SahaSajib2015Eosr}.

Even though standardization has been successfully applied in brain imaging \cite{deGooijer2013_14,Coito2019_18,LiRui2021_7}, its theoretical foundation remains unexplored and thus it is yet unclear when this method performs well and what are the conditions under which we can recover an exact source location. This work aims to provide the theoretical framework for the application of standardization to problems that are modeled by the Poisson differential equation. In particular, we consider the standardization technique within a general discretized Poisson-type forward problem and a general Gaussian process regression model for inversion. In practice, we analyze the forward model-induced localization bias related to the $L2$-norm data fitting and how standardization can reduce it in simple settings. In this study, we assume that there is only one source to be localized and the discretized forward model describes reality reasonably well, i.e., there is no uncertainty related to it. In our analysis in Section \ref{sc:bias_reduction}, it means that the measurements can be written as $\underline{\bf y}=a{\bf L}_{\underline{k}}+\underline{\bf q}$, where ${\bf L}_{\underline{k}}$ is random $\underline{k}$th column of the system matrix and $a$ is source strength. As we will see later, we can bypass the limitation of one source by generalizing the standardization technique for every zero-mean Gaussian prior and using an educated guess about the prior via Kalman filtering.

As our main contributions, we formalize a Bayesian extension of the technique itself, derive the general conditions for perfect localization via the maximum of the reconstruction, and examine the noise robustness of the localization by deriving a lower bound probability for perfect localization for Gaussian distributed measurements. In numerical examples, we demonstrate how multiple simultaneous and correlated sources can be recovered and tracked using the generalized formulation of standardization which is an impossible task with the original standardization \cite{WagnerMichael2004sLORscr,MohdZulki2022sLORmultisource}. We also demonstrate how the lower bound relates to the observed localization results.

The paper is structured as follows. The forward problem and standardized inversion method will be introduced in Section \ref{sc:methods}. The forward model-induced localization bias and how it is reduced by standardization is examined in Section \ref{sc:bias_reduction}, where we use the measurement noise realization ${\bf q}={\bf 0}$ to simplify the terminology and notations. The necessary conditions for perfect localization via generalized standardization are presented in Section \ref{sc:generalization}. The analysis is then further strengthened to consider stochastic Gaussian distributed data in Sections and \ref{sc:noiserobust}, where a lower probability bound for perfect localization is presented. Computational examples of unbiased localization and noise durability are given in Section \ref{sc:ComputationExampl}. Finally, this paper ends with a discussion and conclusion.

\section{Methods}\label{sc:methods}

\subsection{Standardized low-resolution brain tomography as probabilistic model}\label{sc:sLORETA}

In this section, we give a probabilistic interpretation of the standardization technique introduced in \cite{PascualMarqui2002}. Simultaneously, we generalize the original underlying prior model one step further assuming independent but not identically Gaussian distributed sources. Standardization is a post-hoc weighting of minimum norm estimate \cite{HamalainenMNE} developed to overcome its localization bias. The method aims to have the highest magnitude of the reconstruction located where the true source lies, thus producing unbiased estimations when measurement noise is neglectable, assuming no modeling uncertainties, i.e., both observation and Bayesian models are perfect. 

Let us assume that there is a random source spread and placement $\underline{\bf x}\sim \pazocal{N}({\bf 0},\Gamma)$ and noise vector $\underline{\bf q}\sim \pazocal{N}({\bf 0},C)$ forming the randomized measurements $\underline{\bf y}=L\underline{\bf x}+\underline{\bf q}$. The task is to find the most probable originator for the source spread. This corresponds to the minimization problem of the form (\ref{eq:MNE_minimization_problem}) in which $L$ has been replaced with $L\Gamma^{1/2}$ as the forward operator, where $[\Gamma^{1/2}]_{ii}=\Gamma^{1/2}_{ii}$, and identity matrix as the Tikhonov  regularization matrix. Because of this, we are able to rewrite the measurements as $\underline{\bf y}=L\Gamma^{1/2}\underline{\bf u}+\underline{\bf q}$ where naturally $\underline{\bf u}\sim\pazocal{N}({\bf 0},I)$.

A rigorous definition for the post-hoc weights can be given by filling the gaps in the reasoning in \cite{PascualMarqui2002} using probability theory: Since the resolution matrix itself is an intermodel concept, we should also consider the weighted reconstruction as such, and therefore, we could derive the distribution of measurements using marginalized likelihood, i.e.,
\begin{eqnarray}
    \pi(\underline{\bf y})=\int_\Omega \pazocal{N}(L\Gamma^{1/2}{\bf u}\mid C)\pazocal{N}({\bf u}\mid I)\,\mathrm{d}{\bf u}=\pazocal{N}({\bf 0}\mid L\Gamma L^\mathrm{T}+C)
\end{eqnarray}
from which one could define the distribution for the MAP estimator of ${\bf u}$ as
\begin{eqnarray}
    \pi(\underline{\hat{\bf u}})=\pi(\Gamma^{-1/2}L^\dagger\underline{\bf y})=\pazocal{N}\left({\bf 0},\Gamma^{-1/2}R\Gamma^{1/2} \right).
\end{eqnarray}

In addition, we should remark that $\underline{\hat{\bf u}}$ is not a MAP estimator of the original source localization problem. The standardized estimation is then projected back by $\Gamma^{1/2}$. Consequently, by component-wise standardizing $\hat{\bf u}$ for some sampled ${\bf y}$, we obtain
\begin{eqnarray}
    {\bf e}_k^\mathrm{T}\Gamma^{-1/2}{\bf z} = \frac{{\bf e}_k^\mathrm{T}\Gamma^{-1/2}L^\dagger {\bf y}}{\sqrt{{\bf e}_k^\mathrm{T}\Gamma^{-1/2}R\Gamma^{1/2} {\bf e}_k}},
\end{eqnarray}
where ${\bf e}_k$ is $k$th unit basis vector of $\mathbb{R}^n$ and ${\bf z}\in\mathbb{R}^n$ is the parameter vector corresponding to the original source localization problem. Since $\Gamma$ is diagonal,
\begin{eqnarray}
    z_k=\frac{{\bf e}_k^\mathrm{T}\Gamma L^\mathrm{T}\left(L\Gamma L^\mathrm{T}+C\right)^{-1} {\bf y}}{\sqrt{{\bf e}_k^\mathrm{T}R {\bf e}_k}}=R_{kk}^{-1/2}\hat{x}_k
\end{eqnarray}
of which the largest in absolute value is chosen to indicate the estimated position of the source.

The standardization technique is analogous to the Gaussian standardization of the previous MAP estimation applied to each vector component individually. It means that each projection $[L\Gamma^{1/2}]_k$ to measurement space is equally many standard deviations away from the mean, which is zero. This is equivalent to Z-scores \cite{WagnerMichael2004sLORscr,LealAlberto2008sLORscr}, where a set of Gaussian random variables from different scales are unified to ensure a fair comparison. 

As standardization is a weighting of the solution vector depending only on the model parameters, we are able to write the standardized minimum problem as the following minimization problem
\begin{eqnarray}\label{eq:sLORproblem}
    \min_{\bf z\in\mathbb{R}^n}\left\lbrace ({\bf y}-LT{\bf z})^\mathrm{T}C^{-1}({\bf y}-LT{\bf z})+{\bf z}^\mathrm{T}T^{\mathrm{T}}\Gamma^{-1}T{\bf z}\right\rbrace,
\end{eqnarray}
where $T=\mathrm{Diag}\left(R\right)^{1/2}$ and $\mathrm{Diag}(R)$ denotes a diagonal matrix with $\mathrm{diag}(R)$ as its diagonal elements. It should be noted that due to the formulation of the standardization, the technique works best in locating one source. The shortcomings of the method in locating several simultaneous sources have been demonstrated before \cite{MohdZulki2022sLORmultisource}.

\section{Localization bias of BMNE and its reduction}\label{sc:bias_reduction}
\begin{figure}
\flushright
    \begin{minipage}{0.4\textwidth}
        \includegraphics[width=0.95\textwidth]{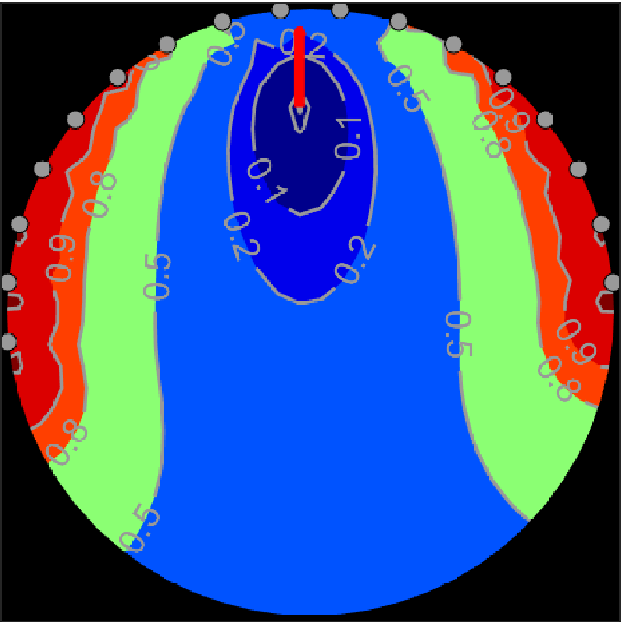}
    \end{minipage}\begin{minipage}{0.4\textwidth}
        \includegraphics[width=0.95\textwidth]{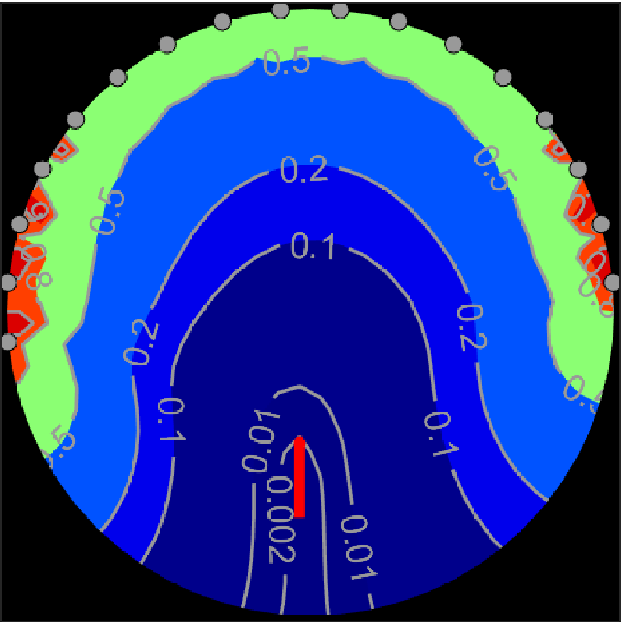}
    \end{minipage}
    \caption{Here we plot the normalized scalar rejection $1-\mathrm{proj}_{\bf y}({\bf L}_k\mid I)$ for two different source configurations. A homogeneous conductivity disk where measurement sensors are on the upper half of the disk (grey balls) is used as the modeling scheme. Low values (dark blue area) are achieved when the columns of the system matrix associated with this area are almost or completely parallel (as vectors) to measurements caused by the source indicated by the red arrow. Values near 1 are achieved when the columns are normal to measurement. The exact value 1 is indicated by dark red coloring. Decomposing each column of the system matrix to a vector that can explain the measurements (that is parallel in the optimal case) and nuisance part (normal) helps us to understand the depth bias obtained from the estimator. The left image shows the value contours for a source near the sensors, and the right-hand-side one represents a situation where the source is far from the sensors.}
    \label{fig:ConduDisk}
\end{figure}

In this section, we examine which model parameters are responsible for the solution bias of BMNE, i.e., which one determines the component of the reconstruction vector with the largest absolute value. For that reason, we need to consider the orthogonal decomposition of system matrix columns to components parallel to measurement $y$, and to those normal to it. We could interpret these components as how much a column ${\bf L}_k$ of the system matrix contributes to explaining the measurements. To find the component of a vector in the direction of the measurement, we use normalized scalar projection (NSP)
\begin{eqnarray}\label{eq:NSP}
    \mathrm{proj}_{\bf y}({\bf L}_k\mid C)=\frac{{\bf y}^\mathrm{T}C^{-1}{\bf L}_k}{\sqrt{{\bf y}^\mathrm{T}C^{-1}{\bf y}{\bf L}_k^\mathrm{T}C^{-1}{\bf L}_k}}
\end{eqnarray}
presented for $k$th system matrix column, from which we can define normalized scalar rejection as $1-\mathrm{proj}_{\bf y}({\bf L}_k\mid C)$. Utilizing the latter, we can obtain that orthogonal sets exist, for example, in a homogenous conductivity disk in Figure \ref{fig:ConduDisk}. From the illustration, we see how the perpendicular set is larger for a source far from the sensors, indicating that the accurate localization of such a source is more difficult than a source near sensors.
\subsection{Example: localization bias of Minimum norm estimate in three source locations}

Let us find the reconstruction $\hat{\bf x}=[\hat{x}_1,\: \hat{x}_2,\: \hat{x}_3]$ of a point source. As an illustration, assume that the system matrix can be written as $L=[{\bf L}_1,\: a{\bf L}_1,\: {\bf L}_3]$, where $a\in\mathbb{R}$. Now, assuming that we have a perfect forward model and the true source is at the location of component 1. To simplify the notation, we define ${\bf y}={\bf L}_1$. In other words, the source is realized at the first node with unit strength, and the measurement noise is realized to zero. We also define Mahalanobis distance as 
\begin{eqnarray}
    \rho({\bf u},{\bf v}\mid \Sigma)=\sqrt{({\bf u}-{\bf v})^\mathrm{T}\Sigma^{-1}({\bf u}-{\bf v})}.
\end{eqnarray}
By utilizing Mahlanobis distance, we can decompose another system matrix entry ${\bf L}_3$ to parallel and orthogonal part with respect to ${\bf L}_1$ such that
\begin{eqnarray}
    \rho({\bf L}_3,{\bf 0}\mid C)^2=\rho(b{\bf L}_1,{\bf 0}\mid C)^2+\rho({\bf L}_1^\perp,{\bf 0}\mid C)^2
\end{eqnarray}
for some $b\in\mathbb{R}$. With this notation, the minimization problem (\ref{eq:MNE_minimization_problem}) can be written as
\begin{eqnarray}
    \fl\min_{\bf x\in\mathbb{R}^3}\left\lbrace \rho({\bf L}_1,{\bf L}_1x_1+a{\bf L}_1x_2+b{\bf L}_1x_3\mid C)^2+\rho({\bf 0},{\bf L}_1^\perp x_3\mid C)^2+\rho({\bf 0},{\bf x}\mid \Gamma)^2\right\rbrace.
\end{eqnarray}
Moreover, when the $3\times 3$ prior covariance is diagonal,
\begin{eqnarray}
    \Gamma = \left[\matrix{
        p_1&0&0\cr 0&p_2&0\cr 0&0&q \cr}\right],
\end{eqnarray}
the solution is of the form
\begin{eqnarray}
    \left[\matrix{
        \hat{x}_1\cr \hat{x}_2 \cr \hat{x}_3\cr}\right]=\left[\matrix{
        p_1{\bf L}_1^\mathrm{T}(L^\mathrm{T}\Gamma L+C)^{-1}{\bf L}_1\cr
        ap_2{\bf L}_1^\mathrm{T}(L^\mathrm{T}\Gamma L+C)^{-1}{\bf L}_1\cr
        b(q^{-1}+({\bf L}_1^\perp)^\mathrm{T}C^{-1}{\bf L}_1^\perp)^{-1}{\bf L}_1^\mathrm{T}(L^\mathrm{T}\Gamma L+C)^{-1}{\bf L}_1\cr}\right].
\end{eqnarray}
We can easily see that the localization estimation is correct, i.e., absolute value $\hat{x}_1$ is the largest if and only if
\begin{eqnarray}\label{eq:BiasIneq}
    p_1>\left|a\right|p_2\quad\textnormal{and}\quad p_1>\left|b\right|(q^{-1}+({\bf L}_1^\perp)^\mathrm{T}C^{-1}{\bf L}_1^\perp)^{-1}.
\end{eqnarray}
Relations between the expressions $p_1$, $\left|a\right|p_2$, and $\left|b\right|(q^{-1}+({\bf L}_1^\perp)^\mathrm{T}C^{-1}{\bf L}_1^\perp)^{-1}$ in the inequalities dictate the biasing of the inverse solution. In effect, the system matrix columns model the measurements produced by sources in different locations. Due to the maximum principle, the columns corresponding to sources near the boundary will have larger absolute values than a column corresponding to a deeper interior location. If the estimated variables are assumed to be identically distributed, meaning that there is a positive number $\delta$ that $\Gamma =\delta I$, then the location of the maximum entry of the reconstruction mostly depends on the system matrix values if the system matrix columns are not distinct enough, i.e., NSP is high between columns of the system matrix.

\subsection{Standardization reduces solution bias}

In this section, we show that the standardization technique introduced in \ref{sc:sLORETA} reduces the bias of BMNE presented in the previous section. In practice, we show that the set of parameter vector components attaining the maximum magnitude is analogous to the set of individual system matrix columns that are parallel with given measurements ${\bf y}$. Here we denote the set of column indices such that the corresponding system matrix columns are parallel to measurement ${\bf y}$ by 
\begin{eqnarray}\label{eq:SolutionSpaceClassic}
    \pazocal{S}({\bf y})=\left\lbrace k\in [1,n]\colon \mathrm{proj}_{\bf y}({\bf L}_k\mid C)=1\right\rbrace.
\end{eqnarray}
In this example, it means that there exists a real number $s_k$ such that ${\bf L}_k=s_k{\bf y}$, when $k\in \pazocal{S}({\bf y})$.
Standardization weights are defined as the reciprocal of the square root of diagonal elements of the resolution matrix 
\begin{eqnarray}
    R_{kk}=\Gamma_{kk}{\bf L}_k^\mathrm{T}(L\Gamma L^\mathrm{T}+C)^{-1}{\bf L}_k.
\end{eqnarray}
If we now assume that the measurement data is produced by a unit-strength source at location $i$ with realization ${\bf y}={\bf L}_i$, we obtain that every standardized component corresponding to the set $\pazocal{S}({\bf y})$ have the form
\begin{eqnarray}
    \left|z_k\right|=\left|\frac{s_k\Gamma_{kk}{\bf L}_i^\mathrm{T}(L\Gamma L^\mathrm{T}+C)^{-1}{\bf L}_i}{\sqrt{s_k^2\Gamma_{kk}{\bf L}_i^\mathrm{T}(L\Gamma L^\mathrm{T}+C)^{-1}{\bf L}_i}}\right|=\sqrt{\Gamma_{kk}{\bf L}_i^\mathrm{T}(L\Gamma L^\mathrm{T}+C)^{-1}{\bf L}_i},
\end{eqnarray}
meaning that the component values are only distinguished by the values of the diagonal prior but not by the magnitude of the corresponding elements of the system matrix. Moreover, the standardized reconstruction has its maximum component at every index $k\in\pazocal{S}({\bf y})$ as shown next.
\begin{lemma}
    The maximum of 
    \begin{eqnarray}
        \max_{{\bf u}\in\mathbb{R}^n}\left|\frac{{\bf u}^\mathrm{T}\Sigma^{-1}{\bf v}}{({\bf u}^\mathrm{T}\Sigma^{-1}{\bf u})^{1/2}}\right|    
    \end{eqnarray}
    is obtained for ${\bf u}$ parallel with {\bf v}, where {\bf v} is non-zero vector and $\Sigma$ is positive definite.
\end{lemma}
\begin{proof}
We can see by straight calculation that
\begin{eqnarray}\label{eq:sup}
    \max_{\bf u}\left|\frac{{\bf u}^\mathrm{T}\Sigma^{-1}{\bf v}}{({\bf u}^\mathrm{T}\Sigma^{-1}{\bf u})^{1/2}}\right|=\max_{\hat{\bf u}}\left|\frac{\hat{\bf u}^\mathrm{T}\hat{\bf v}}{(\hat{\bf u}^\mathrm{T}\hat{\bf u})^{1/2}}\right|=\left\|\hat{\bf v}\right\|_2=\left({\bf v}^\mathrm{T}\Sigma^{-1}{\bf v}\right)^{1/2},
\end{eqnarray}
due to the positive definity of the matrix $\Sigma$, and the maximum is attained when ${\bf u}$ is parallel to ${\bf v}$. 
\end{proof}
It follows that in our case of the Bayesian model, the maximum reconstruction vector component is the one having maximum prior variance $\Gamma_{kk}$, where $k\in\pazocal{S}({\bf y})$. Therefore, the prior model can be used to guide the solution to a certain region based on the prior information on the source location and to avoid the possible case where the solution is non-unique due to $\pazocal{S}({\bf y})$ containing more than one index. The alternative unbiasing technique called {\em exact low-resolution electromagnetic tomography} (eLORETA) \cite{Pascual-Marqui2011_eLORETA} is less flexible in this regard because the weighted regularization is fixed, and thereby, does not allow any alternation of the solution via prior information.

Moreover, we can connect standardization with Mahalanobis distance by noticing
\begin{eqnarray}
    \sqrt{\Gamma_{jj}{\bf L}_j^\mathrm{T}(L\Gamma L^\mathrm{T}+C)^{-1}{\bf L}_j}=\sqrt{\Gamma_{jj}}\rho({\bf 0},{\bf L}_j\mid L\Gamma L^\mathrm{T}+C),
\end{eqnarray}
where the distance is taken with respect to the distribution of measurements as a marginal distribution from the likelihood. Standardization by this Mahalanobis distance can be interpreted so that each projection ${\bf L}_j$ from the source space to the measurement space is equally achievable by the distribution of measurements, i.e., they are equally many standard deviations away from the mean.

\section{Generalized case: Gaussian distributed source configuration}\label{sc:generalization}

To expand the theory presented in Section \ref{sc:sLORETA}, we consider general prior covariance instead of diagonal ones. For that reason, we use the matrix square roots, calculated utilizing Schur factorization \cite{BjorckA1983ASqrtm}, in the modified system matrix $L\Gamma^{1/2}$. One could calculate Cholensky decomposition $\Gamma = MM^T$ and use the left Cholensky factor $M$ instead, but we make our decision for the sake of simplicity in notations and this particular decomposition is suggested to be used \cite{Pascual2007sqrtm} due to its symmetry. In this section, we show that we have a similar index set for maximizing components for reformulated problems as the one presented in the previous section. As in Section \ref{sc:sLORETA}, we consider randomized measurements written as $\underline{\bf y}=L\Gamma^{1/2}\underline{\bf u}+\underline{\bf q}$.

 Let us define the decomposition of the new system matrix first, denoting
\begin{eqnarray}\label{eq:ParallelsubsystemK}
    \hat{K}=\left[L\Gamma^{1/2} \right]_k,\quad k\in\hat{\pazocal{S}}\left({\bf y}\right),
\end{eqnarray}
where
\begin{eqnarray}
    \hat{\pazocal{S}}\left({\bf y}\right)=\left\lbrace k\in [1,n]\colon 
    \mathrm{proj}_{\bf y}\left(\left[ L\Gamma^{1/2}\right]_k\mid C\right)=1\right\rbrace,
\end{eqnarray}
and
\begin{eqnarray}
    \hat{H}=\left[L\Gamma^{1/2}\right]_k,\quad k\in[1,n]\setminus\hat{\pazocal{S}}\left({\bf y}\right).
\end{eqnarray}

The following theorem states that the index set $\hat{\pazocal{S}}({\bf y})$ is, in fact, the set of nodes having maximum reconstruction value for this newly formulated problem.
\begin{theorem}
    Consider a standardized minimum norm problem 
    \begin{eqnarray}
        \min_{{\bf z}\in\mathbb{R}^n}\left\lbrace\rho({\bf y},\hat{K}\hat{T}{\bf z}+\hat{H}\hat{T}{\bf z}\mid C)^2+\rho({\bf 0},\hat{T}{\bf z}\mid I)^2\right\rbrace,
    \end{eqnarray}
    where 
    \begin{eqnarray}
        \hat{T}=\mathrm{Diag}\left(\left[\matrix{
\hat{K}^\mathrm{T} \cr \hat{H}^\mathrm{T}\cr}\right]\left ( \left[\matrix{
\hat{K} & \hat{H} \cr
}\right] \left[\matrix{
\hat{K}^\mathrm{T} \cr \hat{H}^\mathrm{T}
}\right]+C \right )^{-1}\left[\matrix{
\hat{K} & \hat{H} \cr
}\right]\right)^{1/2}.
    \end{eqnarray}
    If there exists a particular solution ${\bf v}_0$ for $\hat{K}{\bf v}_0={\bf y}$, where the sub-system matrix is given by Equation \ref{eq:ParallelsubsystemK}, then the maximum components of the reconstruction vector ${\bf z}$ corresponds to the system matrix columns that can be written as $\hat{K}{\bf v}$ and ${\bf v}$ belongs to the subspace $\mathrm{span}({\bf v_0})$. 
\end{theorem}
\begin{proof}
The problem of finding the maximizing index
\begin{eqnarray}
    \mathrm{arg\: max}_{i\in[1,n]}\left\lbrace\frac{{\bf A}_i^\mathrm{T}\Sigma^{-1}{\bf y}}{\sqrt{{\bf A}_i^\mathrm{T}\Sigma^{-1}{\bf A}_i}}\right\rbrace,
\end{eqnarray}
where ${\bf A}_i$ is the $i$th column vector of the matrix $A=[\hat{K}\: \hat{H}]$ and $\Sigma=A A^\mathrm{T}+C$, can be written as the following maximization problem
\begin{eqnarray}
    \mathrm{maximize\:}{\bf x}^\mathrm{T}A^\mathrm{T}\Sigma^{-1} {\bf y}\quad \mathrm{subject\: to\: }\sqrt{{\bf x}^\mathrm{T}A^\mathrm{T}\Sigma^{-1} A{\bf x}}=c,
\end{eqnarray}
for some $c>0$ and where ${\bf A}_i$ is replaced by $A{\bf x}$. The problem can be solved by first forming the following Lagrangian
\begin{eqnarray}
    \mathcal{L}({\bf x},\lambda)={\bf x}^\mathrm{T}A^\mathrm{T}\Sigma^{-1} {\bf y}-\lambda(\sqrt{{\bf x}^\mathrm{T}A^\mathrm{T}\Sigma^{-1} A{\bf x}}-c)
\end{eqnarray}
and then calculating the derivatives
\begin{eqnarray}
    \frac{\partial \mathcal{L}}{\partial {\bf x}}&=(A^\mathrm{T}\Sigma^{-1} {\bf y})^\mathrm{T}-\frac{\lambda {\bf x}^\mathrm{T}A^\mathrm{T}\Sigma^{-1} A}{\sqrt{{\bf x}^\mathrm{T}A^\mathrm{T}\Sigma^{-1} A{\bf x}}},\\
    \frac{\partial \mathcal{L}}{\partial \lambda}&=c-\sqrt{{\bf x}^\mathrm{T}A^\mathrm{T}\Sigma^{-1} A{\bf x}}.
\end{eqnarray}
By setting these to zero, we get
\begin{eqnarray}
    A^\mathrm{T}\Sigma^{-1} \left({\bf y}-\frac{\lambda}{c}A{\bf x}\right)={\bf 0}.
\end{eqnarray}
However, there is no solution for $k\in[1,n]\setminus\hat{\pazocal{S}}\left({\bf y}\right)$ unless $\Sigma^{-1}({\bf y}-(\lambda/c){\bf A}_k)$ belong to the null-space of $A^\mathrm{T}$ for some $\lambda\in\mathbb{R}$. This again is not possible, since ${\bf y}$ is assumed to be image of $A$. Therefore, there are no local maxima for that index set. Contrarily there exists maximizing index $k\in\hat{\pazocal{S}}\left({\bf y}\right)$ by definition.

In addition, we can write a new maximization problem for the second case, noting that optimum should contain a vector linearly dependent on {\bf y}, and also ${\bf y}^\mathrm{T}\Sigma^{-1} {\bf y}$ does not affect on maximum so that it can be neglected from the problem statement, that is
\begin{eqnarray}
    \mathrm{maximize\:} {\bf x}^\mathrm{T}\hat{H}^\mathrm{T}\Sigma^{-1} {\bf y}\quad \mathrm{subject\: to\: }\sqrt{(\hat{H}{\bf x}+{\bf y})^\mathrm{T}\Sigma^{-1} (\hat{H}{\bf x}+{\bf y})}=c.
\end{eqnarray}
Solving this similarly as above, derivatives of this Lagrangian are
\begin{eqnarray}
    \frac{\partial \mathcal{L}}{\partial {\bf x}}&=(\hat{H}^\mathrm{T}\Sigma^{-1} {\bf y})^\mathrm{T}-\lambda\frac{{\bf x}^\mathrm{T}\hat{H}^\mathrm{T}\Sigma^{-1} \hat{H}+{\bf y}^\mathrm{T}\Sigma^{-1} \hat{H}}{\sqrt{(\hat{H}{\bf x}+{\bf y})^\mathrm{T}\Sigma^{-1} (\hat{H}{\bf x}+{\bf y})}},\\
    \frac{\partial \mathcal{L}}{\partial \lambda}&=c-\sqrt{(\hat{H}{\bf x}+{\bf y})^\mathrm{T}\Sigma^{-1} (\hat{H}{\bf x}+{\bf y})}.
\end{eqnarray}
This yields to the equation
\begin{eqnarray}
    \hat{H}^\mathrm{T}\Sigma^{-1} \left(\left(1-\frac{\lambda}{c}\right){\bf y}-\frac{\lambda}{c} \hat{H}{\bf x}\right)={\bf 0}
\end{eqnarray}
that has the solution $({\bf x},\lambda)=(\mathrm{null}(\hat{H}^\mathrm{T}\Sigma^{-1} \hat{H}),c)=(\mathrm{null}(\hat{H}),c)$ due to rank-nullity theorem, implying that the maximizing system matrix column comes from the vector space $\mathrm{span}({\bf y})$.
\end{proof}
The formulation of the problem of the previous Theorem is analogous to the following, which is written using the original system matrix:
\begin{eqnarray}
        \min_{{\bf z}\in\mathbb{R}^n}\left\lbrace\rho({\bf y},L\hat{T}{\bf z}\mid C)^2+\rho({\bf 0},\hat{T}{\bf z}\mid \Gamma)^2\right\rbrace,
    \end{eqnarray}
    where 
    \begin{eqnarray}
        \hat{T}=\Gamma^{1/2}\mathrm{Diag}\left(\Gamma^{1/2}L^\mathrm{T}(L\Gamma L^\mathrm{T}+C)^{-1}L\Gamma^{1/2}\right)^{1/2}.
    \end{eqnarray}
The proper standardized problem involving the back-projection to the original problem uses the transformation
\begin{eqnarray}\label{eq:GeneralzWeights}
        T=\Gamma^{1/2}\mathrm{Diag}\left(\Gamma^{1/2}L^\mathrm{T}(L\Gamma L^\mathrm{T}+C)^{-1}L\Gamma^{1/2}\right)^{1/2}\Gamma^{-1/2},
    \end{eqnarray}
which yields, in the case of diagonal covariance, the same as the problem \ref{eq:sLORproblem}. Thus, we have presented a natural extension of the technique.

\section{Noise robustness of localization via standardized model}\label{sc:noiserobust}
Here we consider a case where the assumed Gaussian likelihood model holds, i.e., the measurements are Gaussian distributed. We start by assuming the location of the true source is in the $k$th node so that $\underline{\bf y}\sim \pazocal{N}(\hat{\bf L}_k,C)$. For the standardized method to localize the activity correctly, we must have 
\begin{eqnarray}\label{eq:PerfLocIneq}
    \frac{\left|\hat{\bf L}_k^\mathrm{T}\Sigma^{-1}\underline{\bf y}\right|}{\sqrt{\hat{\bf L}_k^\mathrm{T}\Sigma^{-1}\hat{\bf L}_k}}\geq \frac{\left|\hat{\bf L}_i^\mathrm{T}\Sigma^{-1}\underline{\bf y}\right|}{\sqrt{\hat{\bf L}_i^\mathrm{T}\Sigma^{-1}\hat{\bf L}_i}}\quad\mathrm{for\: all\: }i\in [1,n],
\end{eqnarray}
where $\Sigma=\hat{L}\hat{L}^\mathrm{T}+C$. Without loss of generality, we denote the forward model simply by $\hat{L}$ and assuming the prior model as $\pazocal{N}({\bf 0},I)$ but we could either denote the system matrix by $L\Gamma^{1/2}$ a given prior covariance $\Gamma$ and end up with the same result as the previous section implies. The exact probability that this inequality is true can be derived from the difference of generalized chi-distributed random variables. However, the formulation would be really complex and therefore impractical. So instead, we present a lower bound for the probability that the correct element belongs to the set of reconstruction maxima by its index. 

\begin{theorem}\label{thm:NoiseBound}
    Assume Gaussian distributed measurements $\underline{\bf y}\sim \pazocal{N}(\hat{\bf L}_k,C)$ produced by a point source at $k$th location. Let
    \begin{eqnarray}
        \theta\leq 1-\mathrm{proj}_{\hat{\bf L}_k}\left(\hat{\bf L}_i\mid\Sigma\right),
    \end{eqnarray}
    for every $i\in [1,n]\setminus\pazocal{S}(\hat{\bf L}_k)$. Then the following inequality holds for the probability of the $k$th standardized parameter attaining the maximum magnitude
    \begin{eqnarray}
        \fl\mathbb{P}\left(\lbrace k\rbrace\subseteq \mathrm{arg\: max}_{i\in[1,n]}\left\lbrace\frac{\left|\hat{\bf L}_i^\mathrm{T}\Sigma^{-1}\underline{\bf y}\right|}{\sqrt{\hat{\bf L}_i^\mathrm{T}\Sigma^{-1}\hat{\bf L}_i}}\right\rbrace\right)\geq \mathrm{Ga}\left(\theta^2\frac{\lambda_{\mathrm{min}}(\Sigma)\hat{\bf L}_k^\mathrm{T}\Sigma^{-1}\hat{\bf L}_k}{2\lambda_{\mathrm{max}}(C)}\,;\:\frac{m}{2},\:1\right),
    \end{eqnarray}
    where $\lambda_{\mathrm{max}}(\cdot)$ and $\lambda_{\mathrm{min}}(\cdot)$ denotes the largest and smallest eigenvalues of a matrix, respectively, and $\mathrm{Ga}(\cdot;\,\cdot ,\,\cdot)$ is cumulative gamma distribution for the shape parameter $m/2$ and scale parameter 1 with the following expression
\begin{equation}
    \mathrm{Ga}\left(\xi;\frac{m}{2},1\right)=\frac{1}{\Gamma(m/2)}\gamma\left(\frac{m}{2},\xi\right),
\end{equation}
where $\Gamma(\cdot)$ is the gamma function and $\gamma(\cdot,\cdot)$ is the lower incomplete gamma function.
\end{theorem}
\begin{proof}
The proof aims to estimate the probability of 
\begin{eqnarray}
    \frac{\left|\hat{\bf L}_k^\mathrm{T}\Sigma^{-1}\underline{\bf y}\right|}{\sqrt{\hat{\bf L}_k^\mathrm{T}\Sigma^{-1}\hat{\bf L}_k}}\geq \frac{\left|\hat{\bf L}_i^\mathrm{T}\Sigma^{-1}\underline{\bf y}\right|}{\sqrt{\hat{\bf L}_i^\mathrm{T}\Sigma^{-1}\hat{\bf L}_i}}
\end{eqnarray}
for any $i\in[1,n]\setminus\lbrace k\rbrace$.
    Let us denote ${\bf u}=\Sigma^{-1/2}\hat{\bf L}_k/\sqrt{\hat{\bf L}_k^\mathrm{T}\Sigma^{-1/2}\hat{\bf L}_k}$ and $\Sigma^{-1/2}\underline{\bf y}=\Sigma^{-1/2}\hat{\bf L}_k+\underline{\bf f}=\bm{\ell}_k+\underline{\bf f}$ for simplicity. As we can now see that
    \begin{eqnarray}
        \frac{\hat{\bf L}_k^\mathrm{T}\Sigma^{-1}\underline{\bf y}}{\sqrt{\hat{\bf L}_k^\mathrm{T}\Sigma^{-1}\hat{\bf L}_k}}={\bf u}^\mathrm{T}\left(\bm{\ell}_k+\underline{\bf f}\right).
    \end{eqnarray}
    The goal is to estimate an event $Z\geq 0$ for a random number $Z$ defined as
    \begin{eqnarray}
        Z=\left|{\bf u}^\mathrm{T}\left(\bm{\ell}_k+\underline{\bf f}\right)\right|-\left|{\bf v}^\mathrm{T}\left(\bm{\ell}_k+\underline{\bf f}\right)\right|.
    \end{eqnarray}
    By decomposing the normalized system matrix column vector ${\bf v}$ to $\vartheta {\bf u}+\sqrt{1-\vartheta^2}{\bf u}^\perp$, where $\vartheta\in[-1,1]$, ${\bf u}^\perp$ is unit vector and ${\bf u}^\mathrm{T}{\bf u}^\perp=0$, we can obtain a lower bound for the random variable
    \begin{eqnarray}
        Z\geq Z^*(\vartheta)=\left\|\bm{\ell}_k\right\|_2-\left|\vartheta\right|\,\left\|\bm{\ell}_k\right\|_2-\sqrt{1-\vartheta^2}\left\|\underline{\bf f}\right\|_2,
    \end{eqnarray}
    and thus $\mathbb{P}(Z^*(\vartheta)\geq 0)\leq \mathbb{P}(Z\geq 0)$. Considering $Z^*(\vartheta)$ as a function, let us assume that $Z^*(0)$ is positive, i.e., there is more data than noise. Therefore, the expression is positive only when it is bounded by the  zero of $Z^*$:
    \begin{eqnarray}
        \left|\vartheta\right| < \frac{\left\|\bm{\ell}_k\right\|_2^2-\left\|\underline{\bf f}\right\|_2^2}{\left\|\bm{\ell}_k\right\|_2^2+\left\|\underline{\bf f}\right\|_2^2}.
    \end{eqnarray}
    To improve this bound, we obtain
    \begin{eqnarray}
        \frac{\left\|\bm{\ell}_k\right\|_2^2-\left\|\underline{\bf f}\right\|_2^2}{\left\|\bm{\ell}_k\right\|_2^2+\left\|\underline{\bf f}\right\|_2^2}\geq 1-\frac{\left\|\underline{\bf f}\right\|_2}{\left\|\bm{\ell}_k\right\|_2}.
    \end{eqnarray}
    Noting that $\vartheta={\bf u}^\mathrm{T}{\bf v}$, we can write an event in the form
\begin{eqnarray}
    \left\|\underline{\bf f}\right\|_2<\left(1-\left| {\bf u}^\mathrm{T}{\bf v} \right|\right)\left\|\bm{\ell}_k\right\|_2.
\end{eqnarray}
Because we cannot state anything about the vector ${\bf v}$, we set $\theta\in (0,1]$ such that $\theta\leq 1-\left| {\bf u}^\mathrm{T}{\bf v} \right|$ for any ${\bf v}$ not parallel with ${\bf u}$.
Due to the properties of $L2$-norm
\begin{eqnarray}
    \left\|\underline{\bf f}\right\|_2&\leq \left\|\Sigma^{-1/2}C^{1/2}\right\|\cdot \left\|\underline{\hat{\bf f}}\right\|_2\leq\sqrt{\lambda_{\mathrm{max}}\left(\Sigma^{-1}\right)\lambda_{\mathrm{max}}\left(C\right)}\left\|\underline{\hat{\bf f}}\right\|_2\nonumber\\
    &=\sqrt{\frac{\lambda_{\mathrm{max}}(C)}{\lambda_{\mathrm{min}}(\Sigma)}}\left\|\underline{\hat{\bf f}}\right\|_2.
\end{eqnarray}
As the original random number gets smaller values for each sample than the chi-distributed one, the probability of $\left\|\underline{\bf f}\right\|_2<\left(1-\left| {\bf u}^\mathrm{T}{\bf v} \right|\right)\left\|\bm{\ell}_k\right\|_2$ is larger than with the random number at the end of the previous inequality chain, therefore
\begin{eqnarray}
    \mathbb{P}(Z\geq 0)&\geq \mathbb{P}(Z^*(\vartheta)\geq 0)\nonumber \\
    &\geq \mathbb{P}\left(\left\|\underline{\bf f}\right\|_2<\theta\left\|\bm{\ell}_k\right\|_2\right)\geq \mathbb{P}\left(\sqrt{\frac{\lambda_{\mathrm{max}}(C)}{\lambda_{\mathrm{min}}(\Sigma)}}\left\|\underline{\hat{\bf f}}\right\|_2<\theta\left\|\bm{\ell}_k\right\|_2\right).
\end{eqnarray}
The cumulative distribution of the lower bound at the end of the inequality chain can be expressed using gamma distribution with shape parameter $m/2$ and scale parameter 1. So we have
\begin{eqnarray}
    \fl\mathbb{P}\left(\sqrt{\frac{\lambda_{\mathrm{max}}(C)}{\lambda_{\mathrm{min}}(\Sigma)}}\left\|\underline{\hat{\bf f}}\right\|_2<\theta\sqrt{\hat{\bf L}_k^\mathrm{T}\Sigma^{-1}\hat{\bf L}_k}\right)=\mathrm{Ga}\left(\frac{\theta^2\,\hat{\bf L}_k^\mathrm{T}\Sigma^{-1}\hat{\bf L}_k}{2\lambda_{\mathrm{max}}(C)/\lambda_{\mathrm{min}}(\Sigma)};\frac{m}{2},1\right).
\end{eqnarray}
\end{proof}
From the derived inequality, we can see that the lower bound is larger in situations where the noise level is low, and the true source is near the measurement sensors, where the system matrix magnitude is maximized. One could also approximate further by considering a common situation where measurement noise is assumed to be independent and identically distributed with standard deviation $\sigma$. Utilizing the inequality $\hat{\bf L}_k^\mathrm{T}\Sigma^{-1}\hat{\bf L}_k\geq \lambda_\mathrm{min}(\Sigma^{-1})\left\|\hat{\bf L}_k\right\|_2^2=\lambda_\mathrm{max}(\Sigma)^{-1}\left\|\hat{\bf L}_k\right\|_2^2$, one could then write another bound using the signal-to-noise ratio (SNR) defined as the square of the signal mean per variance, i.e., $\left\|\hat{\bf L}_k\right\|_2^2/\sigma^2$, that reads 
\begin{eqnarray}
        \fl\mathbb{P}\left(\lbrace k\rbrace\subseteq \mathrm{arg\: max}_{i\in[1,n]}\left\lbrace\frac{\left|\hat{\bf L}_i^\mathrm{T}\Sigma^{-1}\underline{\bf y}\right|}{\sqrt{\hat{\bf L}_i^\mathrm{T}\Sigma^{-1}\hat{\bf L}_i}}\right\rbrace\right)\geq \mathrm{Ga}\left(\theta^2\frac{\lambda_\mathrm{min}(\Sigma)}{2\lambda_\mathrm{max}(\Sigma)}\cdot\mathrm{SNR}\,;\:\frac{m}{2},\:1\right).
    \end{eqnarray}
Based on this bound, there can be at most 2 \% of noise (at least 35 \unit{\decibel} SNR) in the measurement data to achieve the probability of $0.95$ for the conductivity disk used in Figure \ref{fig:ConduDisk} and independent and identically distributed Gaussian sources. As a drawback of the given bound, without stating something about the rank of the system matrix and decrease of the modeled signal strength along with distance away from the sensor, the bound is generous for far-field sources.

\section{Numerical exmples}\label{sc:ComputationExampl}
This section is divided into three parts. The first part demonstrates the already known fact that the standardized localization is free from the localization bias of BMNE, and the second one shows the advantage of the introduced generalized formulation (Section \ref{sc:generalization}) in multisource localization and tracking. The third section examines the spatial and noise level-related behavior of the lower bound from Theorem \ref{thm:NoiseBound}. In practice, we examine the localization accuracy statistically reflecting the results to the lower bound. Both parts end with their own results sections.

\subsection{Localization of far-field source via BMNE and its standardization}
In this example, we give a demonstration of the source localization in a homogenous conductivity disk. Therefore, this is not reflecting the novelties of this study. Rather, the idea of this example is to familiarize ourselves with the source localization problem at hand and the nature of the BMNE bias.
Let us assume that we have a far-field source, i.e., far away from the sensors that are on the upper half of the disk in Figure \ref{fig:ConduDisk}. The source is located at point (0.2,0.6) and the measurements have 5 \% of Gaussian white noise. The aim is to localize the source via BMNE without and with standardization.

\subsection{results}

\begin{figure}[h!]
    \centering
    \begin{minipage}{0.4\textwidth}
    \begin{center}
        {\bf BMNE}
    \end{center}\vskip0.02cm
        \includegraphics[width=0.95\textwidth]{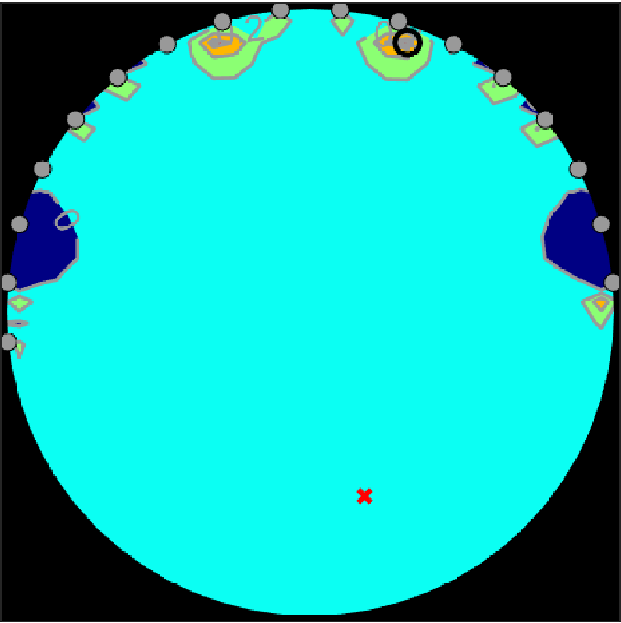}
    \end{minipage}\begin{minipage}{0.4\textwidth}
    \begin{center}
        {\bf sLORETA}
    \end{center}\vskip0.02cm
        \includegraphics[width=0.95\textwidth]{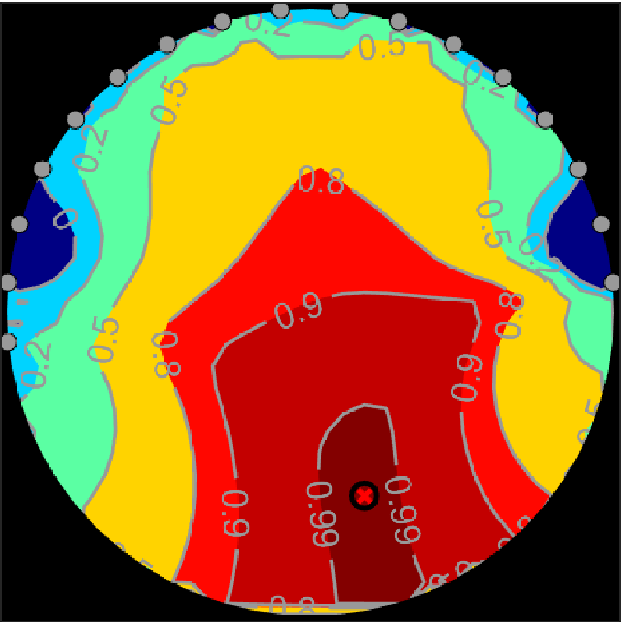}
    \end{minipage}
    \caption{Localization of far-field source located at the position of red cross. The normalized Bayesian minimum norm estimate (BMNE) value map is presented on the left and the normalized value map of standardized BMNE, also called standardized low-resolution brain tomography (sLORETA), is presented on the right. The black ring in both of the pictures shows the location of the maximum reconstruction value where the source is estimated to be. }
    \label{fig:OneSourceLocalization}
\end{figure}

Figure \ref{fig:OneSourceLocalization} presents the heat maps of reconstruction values before and after standardization. As we can see, the BMNE value map is biased toward the sensor locations while after standardization (sLORETA) the localization is perfect. The width of the area of high reconstruction values is notably wider for sLORETA than BMNE.  

\subsection{Unbiasing via standardization and benefits of detailed posterior}
\subsubsection{Simulation setup}
In this section, we aim to track and localize a source time series of two sources via the Kalman filter (KF), sLORETA, and Kalman filter for standardized state variables (SKF) \cite{Lahtinen2024SKF} whose algorithm is presented in \ref{App:SKF}. Here we select our evolution model to be a simple random walk for both Kalman filters.

Let the underlying activity inside the previously used homogenous conductivity disk (Figure \ref{fig:ConduDisk}) obey the following evolution model
\begin{eqnarray}
    {\bf s}_{k+1}=\left[\matrix{
        0.2 & -0.3\cr -0.8&0.3
    \cr}\right]{\bf s}_k+\left[\matrix{
        f((k+1)\Delta t)\cr
        0\cr}\right],
\end{eqnarray}
where $f(t)=\exp(-10^5(t-0.012)^2)\cos(500(t-0.012)+\pi/2)$. The activity is tracked 25 \unit{\milli\second} with sampling rate of 1 \unit{\kilo\hertz}, yielding $\Delta t = 0.001$. The first source (far-field source) is placed at (0,-0.95), and the second (near-field source) is at (-0.4,0.8). For measurement data creation, we add 5 \% of Gaussian noise.
One can clearly see that the activities are not independent, which is why the usage of an advanced posterior model could be appropriate. To estimate the Gaussian posterior along the time series, we apply the Kalman filter to demonstrate the advantage of a detailed posterior model. In the Bayesian sense, the Kalman filter gives a Gaussian posterior for each time step with updated mean and covariance, therefore, based on the section \ref{sc:generalization}, we are able to apply standardization technique to it as a post-hoc weighting, where the weights are updated at every time step.

In the discretization of the analytical forward model, we apply a regular grid of 650 nodes in the disk. For inversion, we form another grid with 30 \% fewer points to avoid an inverse crime. To track the activity, we listen to the location of true sources and track the reconstruction value evolution in those regions.

In localization, we utilize the opposite polarity of the time evolutions. This way, we estimate the location of both sources at every time step by seeking the minimum and maximum value of the reconstruction. Finally, we will assign a Gaussian mixture model to find two clusters that aid in interpreting the localization results.

\subsubsection{Results}
\begin{figure}[h!]
    \flushright
    \begin{minipage}{0.3\textwidth}
    \begin{center}
        {\bf KF}
    \end{center}\vskip0.02cm
        \includegraphics[width=0.9\textwidth]{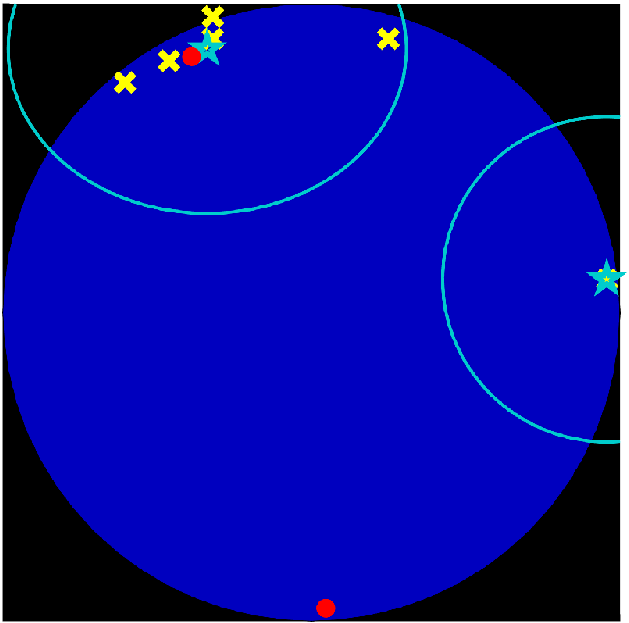}
    \end{minipage}\begin{minipage}{0.3\textwidth}
    \begin{center}
        {\bf SKF}
    \end{center}\vskip0.02cm
        \includegraphics[width=0.9\textwidth]{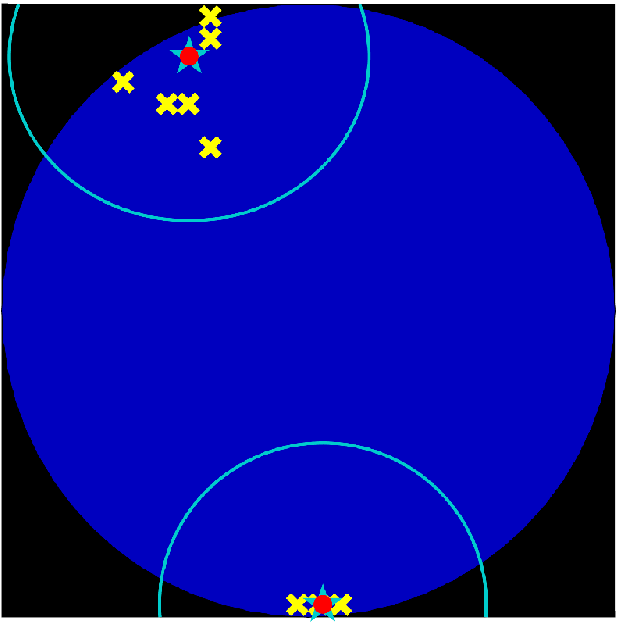}
    \end{minipage}\begin{minipage}{0.3\textwidth}
    \begin{center}
        {\bf sLORETA}
    \end{center}\vskip0.02cm
        \includegraphics[width=0.9\textwidth]{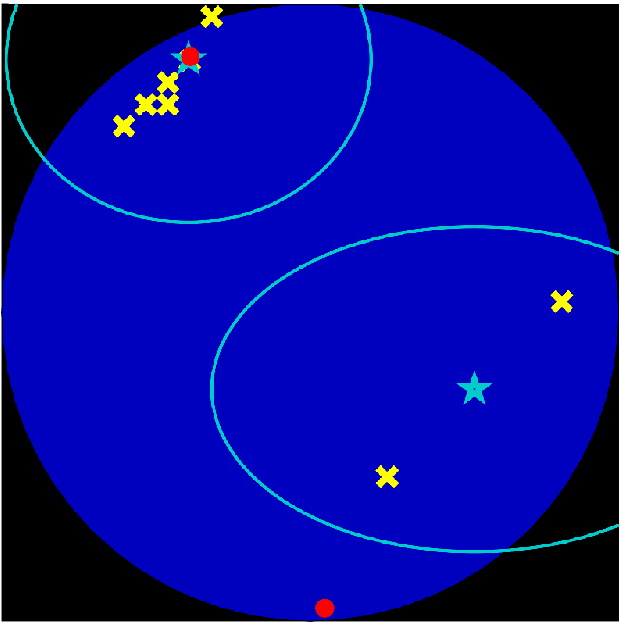}
    \end{minipage}
    \caption{Localization estimates at each time point starting from 6 ms, when the far-field source emerges. Individual estimates are indicated by yellow crosses, red balls show the location of true sources, and a 75 \% confidence ellipse of the 2-component Gaussian mixture model is presented in turquoise. The used methods from left to right are the Kalman filter (KF), standardized Kalman filter (SKF), and standardized low-resolution brain tomography (sLORETA).}
    \label{fig:Localization}
\end{figure}

Considering the localization results in Figure \ref{fig:Localization}, the Kalman filter is able to localize the near-field source with fair accuracy even if the variability in individual localizations at different time steps is large. As KF is analogous to BMNE, the localization bias is prominent since the estimations are accumulated on the upper boundary. sLORETA is able to localize the near-field source with higher accuracy and with less variety than KF. Surprisingly, sLORETA cannot localize the far-field source even if it is detecting signals inside the disk, and therefore, it does not exhibit bias of BMNE. SKF, as a combination of both techniques, is capable of localizing both sources. Variety in near-field localization estimates is slightly lower than with sLORETA as the 75 \% confidence ellipse is a bit smaller. No individual near-field estimation hit the target perfectly but the cluster mean is within sub-mesh accuracy from the true location of the source.

\begin{figure}[h!]
    \centering
    {\bf KF}\\ \vskip0.05cm
    \begin{minipage}{0.55\textwidth}
        \includegraphics[width=0.9\textwidth]{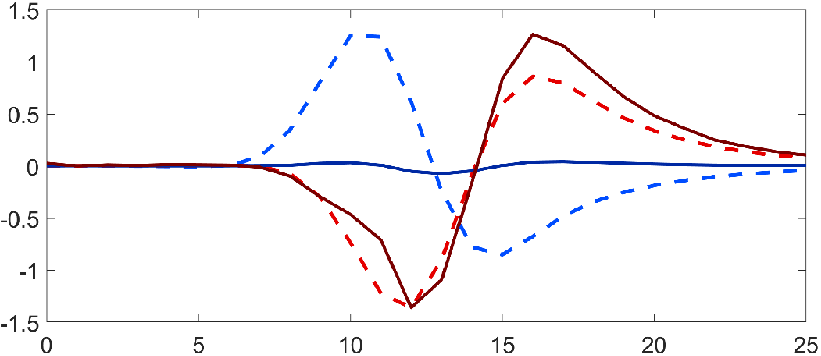}
    \end{minipage}\\ \vskip0.5cm
    {\bf SKF}\\ \vskip0.05cm
    \begin{minipage}{0.55\textwidth}
        \includegraphics[width=0.9\textwidth]{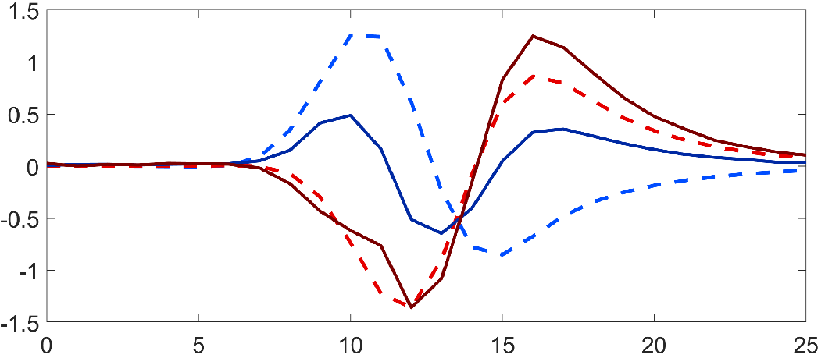}
    \end{minipage}\\ \vskip0.5cm
    {\bf sLORETA}\\ \vskip0.05cm
    \begin{minipage}{0.55\textwidth}
        \includegraphics[width=0.9\textwidth]{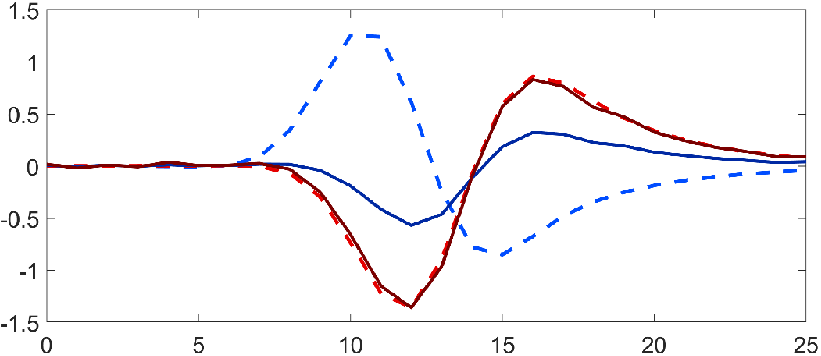}
    \end{minipage}
    \caption{Track of the source activity at the correct location as a function of time that is expressed in milliseconds. Dashed lines represent the true track and solid lines represent the estimated tracks. Blue lines correspond to the far-field source and red lines indicates the near-field source. The used methods from top to bottom are the Kalman filter (KF), standardized Kalman filter (SKF), and standardized low-resolution brain tomography (sLORETA).}
    \label{fig:Tracking}
\end{figure}

\newpage

From the tracking results in Figure \ref{fig:Tracking}, we observe that Kalman filter is tracking the source near sensors, but the relative strength of the far source is suppressed. sLORETA follows the near-field source well, but the far-field source's track follows the other one and seems to not have independent tracks as one would desire. Contrarily, Kalman filter for standardized states can track both of the activities independently at the beginning of the oscillation. The latter part of the far-field track is almost identical to the sLORETA's track. 

\subsection{Localization accuracy from noisy measurements}
\subsubsection{Simulation setup}
The purpose of these experiments is to demonstrate the spatial (I) and SNR (II) dependence of the lower bound from theorem \ref{thm:NoiseBound} and its relation to actual localization accuracy via reconstruction. We have a great interest in how close the bound is to the localization accuracy statistics observed in simulations. 

In example (I), we place a true source in each mesh node once and then calculate 10,000 samples of measurement with Gaussian noise using 5 \% and 15 \% of noise. After that, we use the standardized reconstruction technique to recover the location of the true source while keeping track of how often the maximum absolute value of reconstructions hits the target. From this, we can draw a contour map over the unit disk and compare it with the map of the lower bound.

In example (II), we take three points from the y-axis and calculate the reconstructions and localization success percentage similarly to experiment (I). This time, we use multiple noise levels. The true sources are located from (0,0.90) to (0,0.6) in decrements of 0.1.

\begin{figure}[h!]
    \centering
    \begin{minipage}{0.03\textwidth}
        \rotatebox{90}{5 \%}
    \end{minipage}\begin{minipage}{0.3\textwidth}
    \begin{center}
        {\bf Localizations}
    \end{center}\vskip0.02cm
        \includegraphics[width=0.95\textwidth]{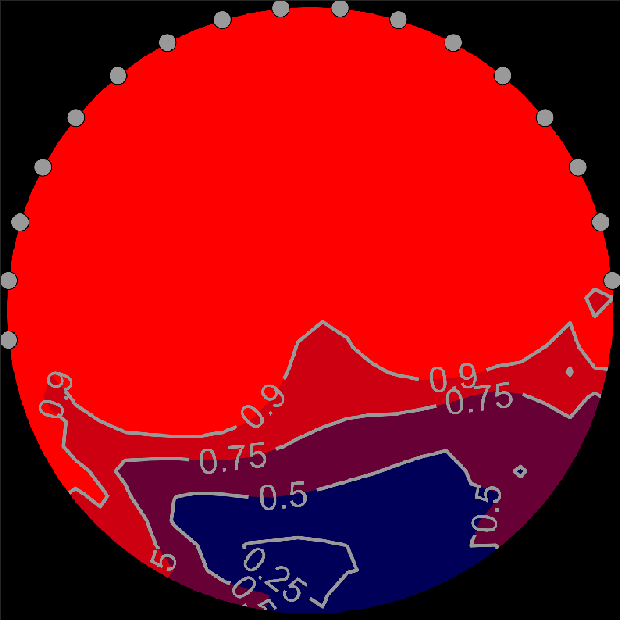}
    \end{minipage}\begin{minipage}{0.3\textwidth}
    \begin{center}
        {\bf Lower bound}
    \end{center}\vskip0.02cm
        \includegraphics[width=0.95\textwidth]{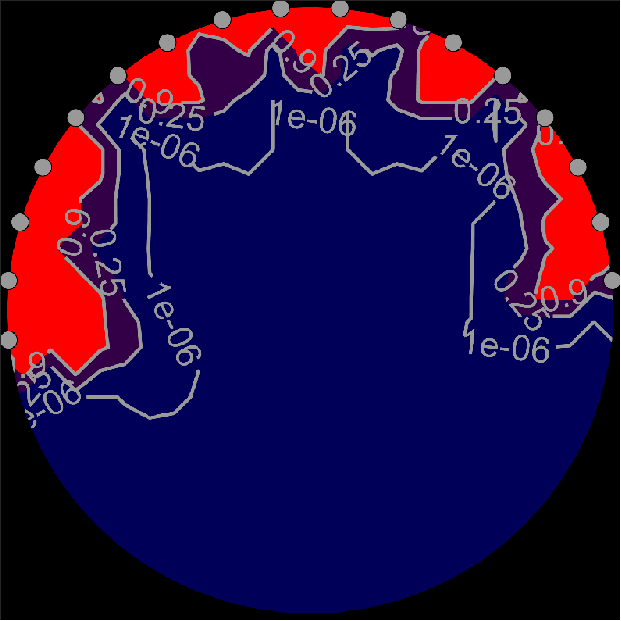}
    \end{minipage}\\
    \hskip0.1cm\begin{minipage}{0.03\textwidth}
        \rotatebox{90}{15 \%}
    \end{minipage}\begin{minipage}{0.3\textwidth}
    \vskip0.02cm
        \includegraphics[width=0.95\textwidth]{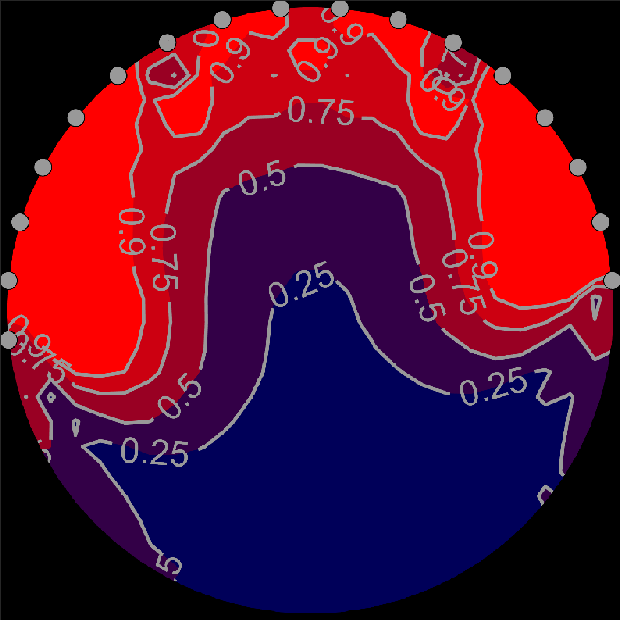}
    \end{minipage}\begin{minipage}{0.3\textwidth}
     \vskip0.02cm
        \includegraphics[width=0.95\textwidth]{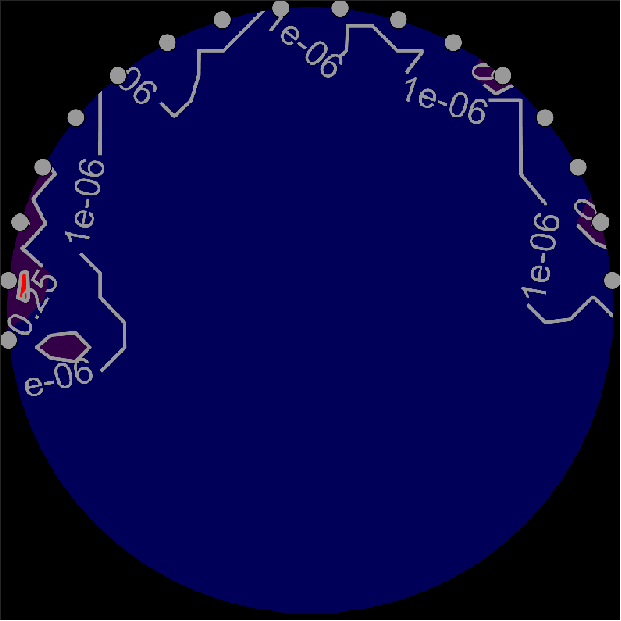}
    \end{minipage}
    \caption{Sampled probability of the standardized inversion to find the correct source location ({\bf left}) mapped over the unit disk and the spatial contours of the lower bound from Theorem \ref{thm:NoiseBound} ({\bf right}). The upper row presents the results when 5 \% of the noise is applied to the measurements and the lower row shows the same results for 15 \% of noise.}
    \label{fig:SpatialProb}
\end{figure}

\subsubsection{Results}
\begin{figure}[h!]
    \centering
   \begin{minipage}{0.25\textwidth}
    \begin{center}
        {\bf (0,0.9)}
    \end{center}\vskip0.02cm
        \includegraphics[width=0.95\textwidth]{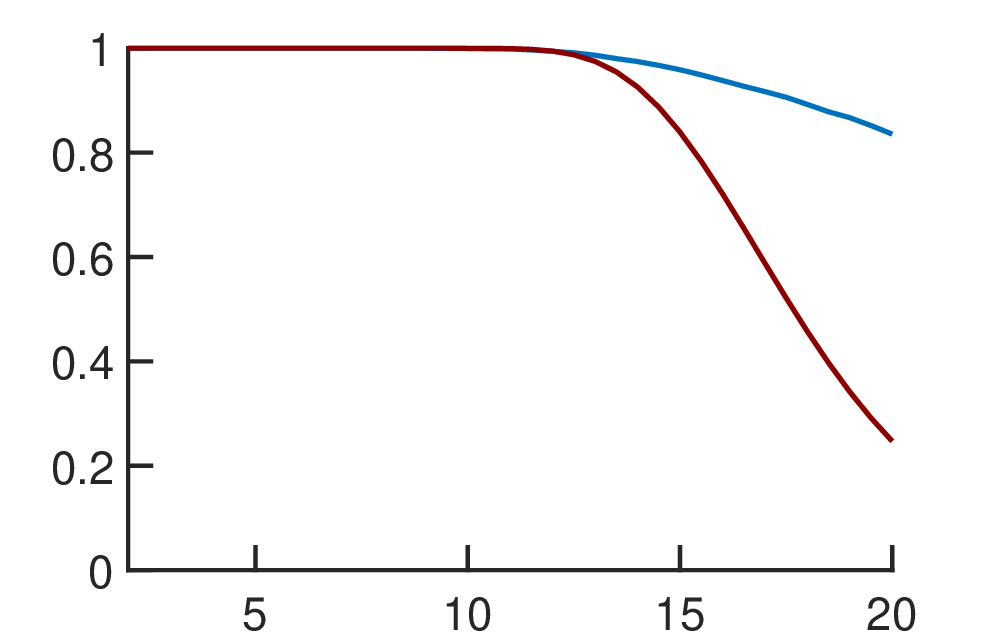}
    \end{minipage}\begin{minipage}{0.25\textwidth}
    \begin{center}
        {\bf (0,0.8)}
    \end{center}\vskip0.02cm
        \includegraphics[width=0.95\textwidth]{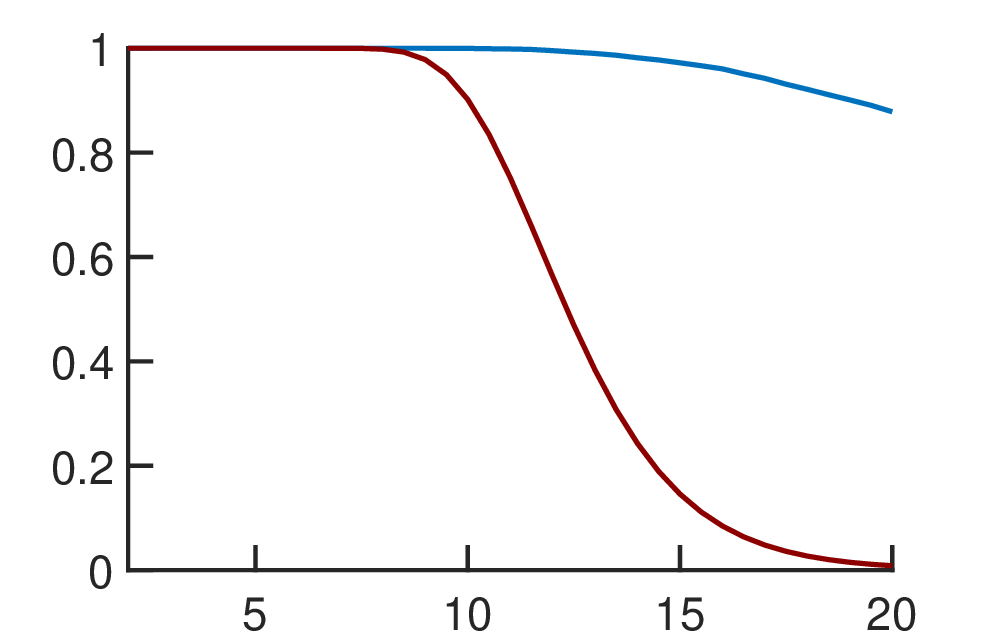}
    \end{minipage}\begin{minipage}{0.25\textwidth}
    \begin{center}
        {\bf (0,0.7)}
    \end{center}\vskip0.02cm
        \includegraphics[width=0.95\textwidth]{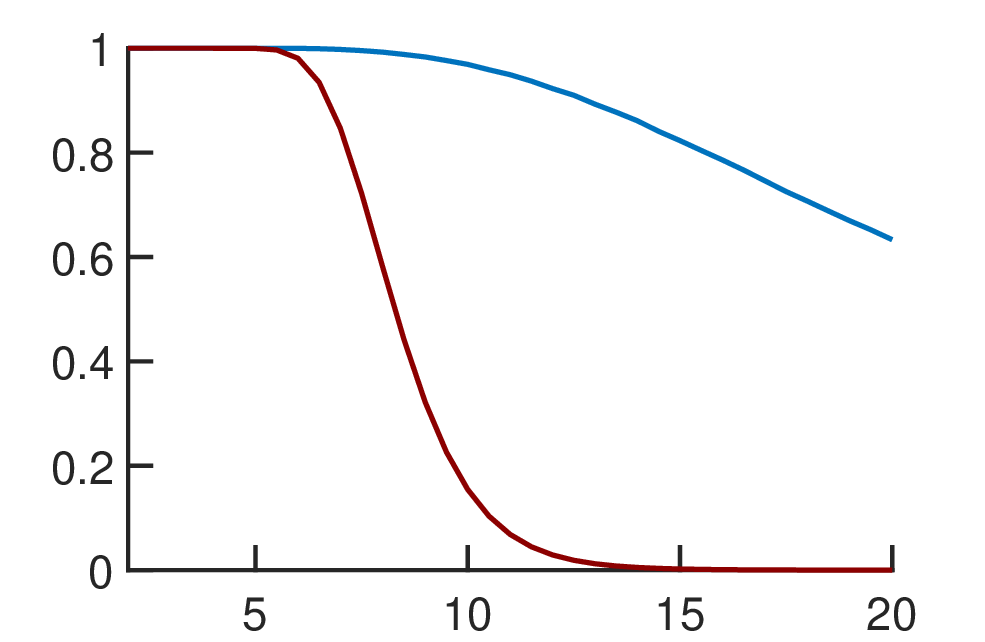}
    \end{minipage}\begin{minipage}{0.25\textwidth}
    \begin{center}
        {\bf (0,0.6)}
    \end{center}\vskip0.02cm
        \includegraphics[width=0.95\textwidth]{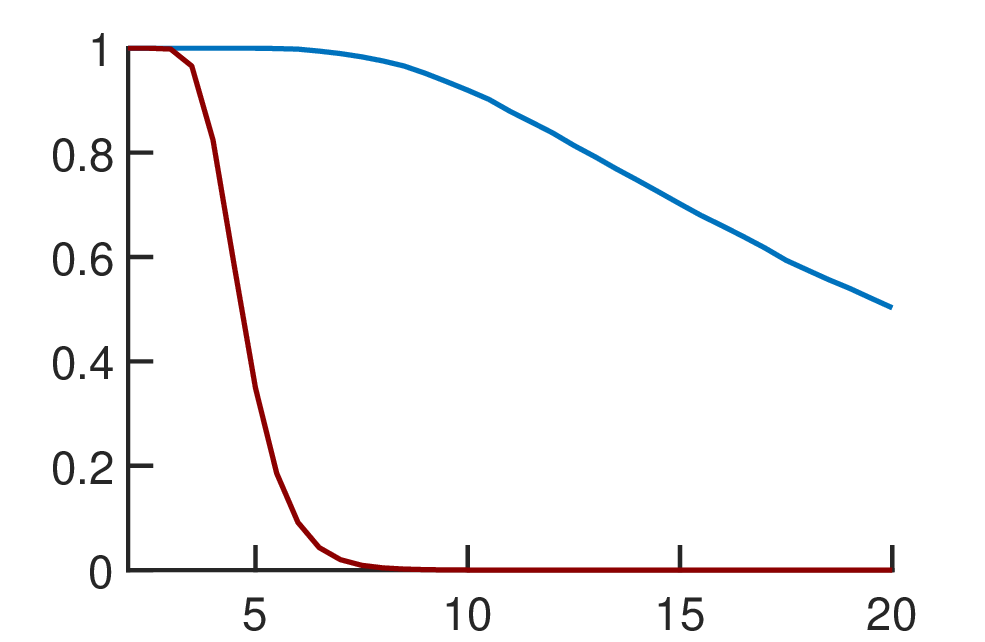}
    \end{minipage}
    \caption{Sampled probability of the standardized inversion to find the correct source location ({\bf blue}) from noise levels from 2 \% to 20 \%. The red curve shows the behavior of the lower bound. Different graphs are calculated at different points of the y-axis and the locations are presented above the graphs}
    \label{fig:SNRProb}
\end{figure}

By looking at the spatial probability maps of the actual reconstruction to localize correctly and the lower bound of that in Figure \ref{fig:SpatialProb}, we can see that the sampled hit rate aligns with the probability bound in the vicinity of the sensors. The lower bound drops quickly along the distance from measurement sensors, much quicker than in the sampled case. In the 5 \% noise case, the standardized reconstruction is accurate almost everywhere with a probability higher than 0.9 except in the region furthest away from the sensors. At 15 \% the region of probability higher than 0.9 probability has shrunk by about half and patterns of lower probability pools emerge near sensors. The lower bound has only a small region of probability higher than 0.9 near the second sensor on the left. The SNR tendency of the sampled probability and the lower bound at different depths in Figure \ref{fig:SNRProb} shows small gaps between probability curves at high SNR (low percent) but the curves diverge faster the deeper the situation is examined. If we say the curves are divergent when the difference is higher than 0.01, the divergence points are at  13, 9, 6, and 3.5 \% for points from 0.9 to 0.6 in decrements of 0.1.

\section{Discussion and conclusions}\label{sc:conclusions}
This paper considers the standardization technique that is well-merited in the neuroelectromagnetic inverse source problem \cite{deGooijer2013_14,Coito2019_18,LiRui2021_7,vandeVelden2023} described here for the general electromagnetic inverse source problem stemming from linearized Poisson problem with Neumann boundary condition. The bias of minimum norm estimation due to the maximum principle is demonstrated as well as how standardization can reduce this bias in the statistical sense by utilizing the orthogonal decomposition of the system matrix. Then, the inverse technique is generalized for any Bayesian model of Gaussian posterior, i.e., Gaussian process regression problem, and its analogy to the original presentation (sLORETA) that assumes independently and identically distributed Gaussian recovering variables. We also derived a lower bound for the probability that the method localizes the source perfectly under Gaussian distributed measurements. The results show that even if the correct source location is found with a forward model that fully describes reality, the set of solutions may be larger than one point if the simulated measurements produced from different sources are not distinct enough or the measurements are significantly noisy. This could, for example, be realized in a situation where the source configuration can have sources both near and far with respect to the measurement sensors, the only difference of which is the strength of the measurement signal. 

Our approach shares similarities with the consideration of Candès {\em et al.} \cite{CandesEmmanuel2006} about perfect reconstructions of a sparse time-discrete signal on the prime-cyclic group $\mathbb{Z}/p\mathbb{Z}$ based on limited information about its frequency components. In this inversion application, the Fourier transform acts as a system matrix. Although the Fourier transform has no near- or far-field so that an inversion would exhibit depth bias, the projections of its column vectors to signal coefficients play a significant role: It is easy to see that in the single signal spike case, the absolute value of normalized scalar projection (NSP) must be equal to one for the correct solution, as for the feasible set $\pazocal{S}({\bf y})$ of sLORETA. The optimization methods studied in the article, i.e., minimization of the L1 norm of the reconstruction and minimization of the total variation of the reconstruction, can be seen as ways to maximize the signal projection. Especially $L1$-optimization penalizes heavily when the projection is not perfect.

Because the technique standardizes the possible solutions in relation to the location, the strength of the signal no longer has an effect. Therefore, one must carefully design the Gaussian prior model to pick the most probable source location from the feasible set provided by standardization. Although it is demonstrated previously that the standardized method utilizing Gaussian posterior cannot recover multiple simultaneous sources when independent and identically distributed sources are assumed as {\em a prior} \cite{WagnerMichael2004sLORscr,MohdZulki2022sLORmultisource}, our results suggest that this is not the case with more complex Gaussian priors. Here we gave the standardization of the Kalman filter as a practical example, which Bayesian generalization has made possible.

Also other valid approaches exist for estimating possibly full prior covariance. For example, if enough data is available, the covariance matrix could be built using cross-validation, where data is split into two groups. One group is used in the model built and the other for validation. An estimation of the prior covariance matrix could also be computed using the maximum likelihood estimation, introducing hyperparameters and hyperprior distributions \cite{ChenLjung2013}, using generalized singular value decomposition (GSVD) \cite{HANSENP1989GSVD}, unbiased predictive risk estimator \cite{Vogel2002UPRE}, or autoregressive models with exogenous input (ARX) \cite{Ljung2021}. All of these could be used with the generalized standardization technique.

Moreover, using the Bayesian extension of standardization, the technique could be coupled with any conditionally Gaussian model having hyperpriors, for example, gamma, inverse gamma \cite{Calvetti2009}, or half-Cauchy distribution \cite{DongYiqiu2023_hyperpriorHorseshoe} as the hyperprior model.

The probability limit for perfect localization explains why high noise robustness is obtained when the true source is located near measurement sensors \cite{DümpelmannMatthias2012sard,SahaSajib2015Eosr}. Even though the presented lower limit compares weakly to the hit rate of practical localization for far-field sources and with particularly high noise levels, the limit could have its uses in the estimation of high p-values for common low noise and near-field cases, where sources of error are mostly related to modeling. The noise level is not necessarily very high in actual measurement procedures due to the averaging over multiple trials.  Moreover, the modeling-based accuracy of localization is overshadowed by practical sources of errors, such as unfavorable sensor positioning and the uncertainty of forward model parameters.

\section*{Acknowledgment}
The author would like to thank his supervisors Prof. Sampsa Pursiainen and Dr. Alexandra Koulouri for discussions on this work. 
The author is supported by the Jenny and Antti Wihuri Foundation.

\appendix
\section{Standardized Kalman Filter}\label{App:SKF}
Kalman filter is a Bayesian method developed to solve the discrete-time filtering problem \cite{Kalman1961}. The solution provided by the Kalman filter is optimal in the sense that it minimizes estimated error covariance. In recent years, Kalman Filter has sparked many applications including, e.g.,  navigation \cite{RaitoharjuMatti2015KalmanPositioning}, robotics \cite{ThrunKalmaninRobotics}, the neuroelectromagnetic inverse source problem when applied for medical EEG data \cite{GalkaAndreas2004KalmanEEG}, and other parameter estimation and tracking in time \cite{XiongRui2014KalmanParametersEst}.

Mathematically, the discrete-time filtering problem is about estimating hidden-states $\lbrace{\bf x}_t\rbrace_{t=1}^T$ using the stochastic observations $\lbrace{\bf y}_t\rbrace_{t=1}^T:={\bf y}_{1:T}$ and linear model connecting states to the observations. In addition, it is assumed that there is a linear Markov model describing the time evolution of the states. If the time-evolution and observation model are Gaussians, i.e., particularly

\begin{eqnarray}
    \underline{\bf y}_t&=L\underline{\bf x}_{t}+\underline{\bf r}_t,\\
    \underline{\bf x}_{t}&=A_t\underline{\bf x}_{t-1}+\underline{\bf q}_t,
\end{eqnarray}
where ${\bf r}_t\sim \pazocal{N}({\bf 0},R_t)$ and ${\bf q}_t\sim \pazocal{N}({\bf 0},Q_t)$.
Then, the mean ${\bf x}_{t\mid t}$ and covariance $P_{t\mid t}$ of the Gaussian posterior $p(\underline{\bf x}_t\mid {\bf y}_{1:t})=\pazocal{N}({\bf x}_{t\mid t},P_{t\mid t})$ can be obtained recursively as follows \cite{SarkkaSimo2013}:

The prediction step, which gives the sufficient statistics of the predictive distribution $p(\underline{\bf x}_t\mid {\bf y}_{1:t-1})=\pazocal{N}({\bf x}_{t\mid t-1},P_{t\mid t-1})$ for time step $t$, is calculated by
\begin{eqnarray*}
    {\bf x}_{t\mid t-1}&=A_t{\bf x}_{t-1\mid t-1}\\
    P_{t\mid t-1}&=A_tP_{t-1\mid t-1}A_t^\mathrm{T}+Q_t
\end{eqnarray*}
and the update step is 
\begin{eqnarray*}
    S_t&=LP_{t\mid t-1}L^\mathrm{T}+R_t\\
    K_t&=P_{t\mid t-1}L^\mathrm{T}S_t^{-1}\\
    {\bf x}_{t\mid t}&={\bf x}_{t\mid t-1}+K_t\left({\bf y}_t-L{\bf x}_{t\mid t-1}\right)\\
    P_{t\mid t}&=P_{t\mid t-1}-K_tS_tK_t^\mathrm{T}\\
    W_t &=P_{t\mid t-1}^{1/2}\mathrm{Diag}\left(P_{t\mid t-1}^{-1/2}K_tS_tK_t^\mathrm{T}P_{t\mid t-1}^{-1/2}\right)^{-1/2}P_{t\mid t-1}^{-1/2}\\
    {\bf z}_{t\mid t}&=W_t{\bf x}_{t\mid t}.
\end{eqnarray*}
The second last step introduces the time-dependent post-hoc weights consisting of inversion of matrix square root of the prior covariance matrix and the diagonal matrix, where $\mathrm{Diag}(\cdot)$ forms a diagonal matrix with the same diagonal elements as in the input matrix. The formulation of the weighting matrix is derived from the equation (\ref{eq:GeneralzWeights}). Vector ${\bf z}_{t\mid t}$ gives the standardized Kalman estimation at time step $t$. For further reading of the application of standardized Kalman filtering, see \cite{Lahtinen2024SKF}.

\section*{References}

\bibliography{sample}

\end{document}